\documentclass[a4paper,12pt]{article}
\usepackage{amsmath}
\usepackage{color}
\usepackage{amssymb}
\usepackage{latexsym}
\usepackage{url}
\numberwithin{equation}{section}

\def\R{{\bf R}}

\def\N{{\bf N}}

\def\d{\displaystyle}
\def\e{{\varepsilon}}

\newtheorem{thm}{Theorem}[section]

\newtheorem{lem}{Lemma}[section]

\newtheorem{rem}{Remark}[section]
\newtheorem{Def}{Definition}[section]

\title{Nonexistence of global solutions of wave equations
with weak time-dependent damping
and combined nonlinearity}
\author{
Ning-An Lai
\footnote{Institute of Nonlinear Analysis and Department of Mathematics,
Lishui University,
Lishui 323000,
China.
}
\footnote{
School of Mathematical Sciences, Fudan University, Shanghai 200433, China.
e-mail: hyayue@gmail.com.}
\quad
Hiroyuki Takamura
\footnote{Department of Complex and Intelligent Systems,
Faculty of Systems Information Science,
Future University Hakodate,
116-2 Kamedanakano-cho,
Hakodate, Hokkaido 041-8655, Japan.
e-mail: takamura@fun.ac.jp.}
}
\date{
\[
\begin{array}{ll}
\mbox{\footnotesize{\bf Keywords:}}
& \mbox{\footnotesize scattering damping, combined nonlinearty, blow-up, lifespan}\\
\mbox{\footnotesize{\bf MSC2010:}}
& \mbox{\footnotesize primary 35L71, secondary 35B44}\\
\end{array}
\]
}

\begin{document}
\maketitle
\begin{abstract}
In our previous two works, we studied the blow-up and lifespan estimates
for damped wave equations with a power nonlinearity of the solution or its derivative,
with scattering damping independently.
In this work, we are devoted to establishing a similar result
for a combined nonlinearity.
Comparing to the result of wave equation without damping,
one can say that the scattering damping has no influence.
\end{abstract}


\section{Introduction}
\par\quad
Recently, the small data Cauchy problem of damped semilinear wave equations with time dependent variable coefficients attracts more and more attention. The works of Wirth \cite{Wir1, Wir2, Wir3} showed that the behavior of the solution of the following linear problem
\begin{equation}
\nonumber
\left\{
\begin{array}{l}
\d u^0_{tt}-\Delta u^0+\frac{\mu}{(1+t)^\beta}u^0_t=0,
\quad \mbox{in}\ \R^n \times[0,\infty),\\
u^0(x,0)=u_1(x),\ u^0_t(x,0)=u_2(x), \quad x\in\R^n,
\end{array}
\right.
\end{equation}
heavily relies on the decay rate $\beta$ and the size of the positive constant $\mu$. Then people get interested in the corresponding nonlinear problem, i.e. the following small data Cauchy problem
\begin{equation}
\label{nonlinear1}
\left\{
\begin{array}{l}
\d u_{tt}-\Delta u+\frac{\mu}{(1+t)^\beta}u_t=|u|^p
\quad \mbox{in}\ \R^n\times[0,\infty),\\
u(x,0)=\e f(x),\ u_t(x,0)=\e g(x), \quad  x\in\R^n,
\end{array}
\right.
\end{equation}
where $\mu>0,\ n\in\N$ and $\beta\in\R$, and $\e$ measures the smallness of the data. Before going on, it is necessary to mention two corresponding nonlinear problems without damping
\begin{equation}
\label{heat}
\left\{
\begin{array}{l}
\d u_{t}-\Delta u=|u|^p
\quad \mbox{in}\ \R^n\times[0,\infty),\\
u(x,0)=\e f(x), \quad  x\in\R^n,
\end{array}
\right.
\end{equation}
and
\begin{equation}
\label{wave}
\left\{
\begin{array}{l}
\d u_{tt}-\Delta u=|u|^p
\quad \mbox{in}\ \R^n\times[0,\infty),\\
u(x,0)=\e f(x),\ u_t(x,0)=\e g(x), \quad  x\in\R^n.
\end{array}
\right.
\end{equation}
For Cauchy problem \eqref{heat} we know that it admits the critical value of $p$
by
\[
p_F(n):=1+\frac 2n,
\]
which is so-called Fujita exponent,
while the one for problem \eqref{wave}
is so-called Strauss exponent $p_S(n)$, which is the positive
root of the quadratic equation,
\[
\gamma(p, n):=2+(n+1)p-(n-1)p^2=0.
\]
\begin{rem}
\lq\lq critical" here means the borderline which divides the domain of $p$ into the blow-up part and the global existence part of the solution.
\end{rem}
\begin{rem}
It is easy to prove that
\[
p_F(n)<p_S(n)\quad\mbox{for}\ n\ge2.
\]
\end{rem}

Now we come back to Cauchy problem \eqref{nonlinear1}. It is interesting to study the relation of the critical exponents among \eqref{nonlinear1}, \eqref{heat} and \eqref{wave}. For $\beta\in [-1, 1)$, due to the works \cite{DLR13, LNZ12, WY17, IO, FIW, II}, we know that it admits the same critical exponent as that of problem \eqref{heat}. For $\beta>1$, since the authors showed blow-up result for $1<p<p_S(n)$ in \cite{LT18}, we may believe that it has the same critical exponent as that of \eqref{wave}.

If we consider the case $\beta=1$ for Cauchy problem \eqref{nonlinear1}, the size of the positive constant $\mu$ should also be taken into account. Generally speaking, if $\mu$ is relatively large, the term $\{\mu/(1+t)\}u_t$ in the equation will have the main influence on the behavior of the solution, which means that this case has the same critical exponent as that of problem \eqref{heat}. See the works \cite{DABI, DL1}. But, if $\mu$ is relatively small, we may conjecture that the influence of $u_{tt}$ will dominate over $\{\mu/(1+t)\}u_t$, which means that the critical exponent is related to $p_S(n)$. See the work \cite{LTW} by the authors and Wakasa for $0<\mu<(n^2+n+2)/\{2(n+2)\}$, which was extended to $0<\mu<(n^2+n+2)/(n+2)$ by Ikeda and Sobajima \cite{IS} and Tu and Lin \cite{TL1, TL2}. Unfortunately, till now we are not clear of the boardline of $\mu$, which determines that the critical power of Cauchy problem \eqref{nonlinear1} with $\beta=1$ will be Fujita or Strauss. We refer the reader to a very recent work by Palmieri and Reissig \cite{P-R}.

In a recent work \cite{LT17} by the authors, we study the blow-up for the small data Cauchy problem
\begin{equation}
\label{nonlinear2}
\left\{
\begin{array}{l}
\d u_{tt}-\Delta u+\frac{\mu}{(1+t)^\beta}u_t=|u_t|^p
\quad \mbox{in}\ \R^n\times[0,\infty),\\
u(x,0)=\e f(x),\ u_t(x,0)=\e g(x), \quad  x\in\R^n.
\end{array}
\right.
\end{equation}
If $\beta>1$, then we showed that the problem has no global solution for $1<p\leq p_G(n)$, where
\[
p_G(n):=\frac{n+1}{n-1},
\]
which denotes the critical exponent for Glassey conjecture. In this work, we are devoted to studying the small data Cauchy problem
with combined nonlinear terms, that is:
\begin{equation}
\label{nonlinear3}
\left\{
\begin{array}{l}
\d u_{tt}-\Delta u+\frac{\mu}{(1+t)^\beta}u_t=|u_t|^p+|u|^q
\quad \mbox{in}\ \R^n\times[0,\infty),\\
u(x,0)=\e f(x),\ u_t(x,0)=\e g(x), \quad  x\in\R^n,
\end{array}
\right.
\end{equation}
where $\beta>1$. Inspired by the work \cite{HZ}, in which Han and Zhou studied the Cauchy problem \eqref{nonlinear3} without damping and obtained the blow-up result for
\begin{equation}\label{HZcond}
\max\left(1, \frac{2}{n-1}\right)<p\leq \frac{2n}{n-1}
\end{equation}
and
\[
1<q<\min\left\{1+\frac{4}{(n-1)p-2}, \frac{2n}{n-2}\right\},
\]
we want to show that whether we have the same blow-up result for Cauchy problem \eqref{nonlinear3}. The difficulty comes from the damping term,
which prevents us from getting the lower bound of some functional by using the test function method, and we overcome this by using a multiplier which was first introduced in the authors \cite{LT18}. Also, due to the damping term, we can't get the blow-up result and lifespan estimate by using Kato's Lemma, and we do it by using an iteration argument similar to that in \cite{LT18}.
\begin{rem}
Hidano, Wang and Yokoyama \cite{HWY} established global existence result for Cauchy problem \eqref{nonlinear3} without damping for $n=2, 3$ and
\[
p>p_G(n), q>q_S(n)~ and~(q-1)\left((n-1)p-2\right)\geq 4.
\]
In the following we are going to find out that whether the global existence result holds for Cauchy problem \eqref{nonlinear3}.
\end{rem}

\section{Main Result}
\par
First we introduce the definition of the solution as follows.
\begin{Def}\label{def1}
 As in \cite{LT18}, we say that $u$ is an energy solution of
 \eqref{nonlinear3} on $[0,T)$
if
\[
u\in \bigcap_{i=0}^1 C^i([0,T),H^{1-i}(\R^n))
\cap C^1((0,T),L^p(\R^n))\cap L^q_{loc}(\R^n\times(0,T))
\]
satisfies $u(x,0)=\e f(x)$ in $H^1(\R^n)$ and
\begin{equation}
\label{energysol}
\begin{array}{l}
\d\quad\int_{\R^n}u_t(x,t)\phi(x,t)dx-\int_{\R^n}\e g(x)\phi(x,0)dx\\
\d\quad+\int_0^tds\int_{\R^n}\left\{-u_t(x,s)\phi_t(x,s)+\nabla u(x,s)\cdot\nabla\phi(x,s)\right\}dx\\
\d\quad+\int_0^tds\int_{\R^n}\frac{\mu u_t(x,s)}{(1+s)^{\beta}}\phi(x,s)dx\\
\d=\int_0^tds\int_{\R^n}|u_t(x,s)|^p\phi(x,s)dx+\int_0^tds\int_{\R^n}|u(x,s)|^q\phi(x,s)dx
\end{array}
\end{equation}
with any $\phi\in C_0^{\infty}(\R^n\times[0,T))$ and any $t\in[0,T)$.
\end{Def}

Employing the integration by parts in \eqref{energysol}
and letting $t\rightarrow T$, we get the weak solution of \eqref{nonlinear3}
\[
\begin{array}{l}
\d\quad\int_{\R^n\times[0,T)}
u(x,s)\left\{\phi_{tt}(x,s)-\Delta \phi(x,s)
-\left(\frac{\mu\phi(x,s)}{(1+s)^{\beta}}\right)_s \right\}dxds\\
\d=\int_{\R^n}\mu\e f(x)\phi(x,0)dx-\int_{\R^n}\e f(x)\phi_t(x,0)dx\\
\d\quad+\int_{\R^n}\e g(x)\phi(x,0)dx+\int_{\R^n\times[0,T)}|u_t(x,s)|^p\phi(x,s)dxds\\
\d\quad+\int_{\R^n\times[0,T)}|u(x,s)|^q\phi(x,s)dxds.
\end{array}
\]
\par
Our main theorem is the following.
\begin{thm}
\label{blowup}
Let $\mu>0$, $\beta>1$ and $n\geq 1$.
Assume that both $f\in H^1(\R^n)$ and $g\in L^2(\R^n)$ are non-negative, compactly supported,
and $g$ does not vanish identically.
Suppose that an energy solution $u$ of \eqref{nonlinear3} on $[0,T)$ satisfies
\begin{equation}
\label{support}
\mbox{\rm supp}\ u\ \subset\{(x,t)\in\R^n\times[0,\infty)\ :\ |x|\le t+R\}
\end{equation}
with some $R\ge1$. If
\begin{equation}
\label{cond1}
p>1
\end{equation}
and
\begin{equation}
\label{cond2}
\left\{
\begin{array}{ll}
1<q<\min\left\{1+\d\frac{4}{(n-1)p-2}, \frac{2n}{n-2}\right\}
&\mbox{for}\quad n\geq 2,\\
1<q & \mbox{for}\quad n=1,
\end{array}
\right.
\end{equation}
then there exists a constant $\e_0=\e_0(f,g,n,p,\mu, \beta, R)>0$
such that $T$ has to satisfy
\begin{equation}
\label{lifespan1}
T\leq C\e^{-2p(q-1)/\{2q+2-(n-1)p(q-1)\}}
\end{equation}
for $0<\e\leq \e_0$, where $C$ is a positive constant independent of $\e$.
\end{thm}

\begin{rem}
We have less restriction for $p$, by comparing the conditions \eqref{cond1} and \eqref{HZcond}, since we use an iteration argument instead of
Kato's type lemma. Which means that we may get blow-up result even for large $p$ but small $q$. What is more, for relatively large $p$ and small
$q$, we can establish an improved lifespan estimate. See Theorem \ref{improvedlifespan} below.
\end{rem}
\begin{rem}
The restriction $q<2n/(n-2)$ for $n\geq 2$ is necessary to guarantee the integrability of the nonlinear term $|u|^q$.
\end{rem}
\begin{rem}
As in \cite{HZ}, we should point out that there exist pairs of $(p, q)$ satisfying
\[
p>p_G(n),\quad q>q_S(n),
\]
but still blow-up will occur. For example, since
\[
\gamma\left(n, 1+\frac{4}{n-1}\right)=-\frac{8}{n-1}<0,
\]
we may choose such an appropriate
pair $(p_0,q_0)$ by setting small constants ,
$\delta_1$ and $\delta_2$, such that
\[
p_0:=\frac{n+1}{n-1}+\delta_1>p_G(n)
\]
and
\[
\begin{aligned}
q_S(n)&<q_0:=1+\frac{4}{n-1+(n+1)\delta_1}<1-\delta_2+\frac{4}{n-1}\\
&<1+\frac{4}{n-1+(n-1)\delta_1}=1+\frac{4}{(n-1)p_0-2}.
\end{aligned}
\]
\end{rem}

\par
We also have an improvement on the estimate of the lifespan
for relatively large $p$ and small $q$ as follows.

\begin{thm}
\label{improvedlifespan}
Let $\mu>0$, $\beta>1$ and $n\geq 2$.
Assume that both $f\in H^1(\R^n)$ and $g\in L^2(\R^n)$ are non-negative, compactly supported,
and $g$ does not vanish identically.
Suppose that an energy solution $u$ of \eqref{nonlinear3} on $[0,T)$ satisfies
\begin{equation}
\label{support2}
\mbox{\rm supp}\ u\ \subset\{(x,t)\in\R^n\times[0,\infty)\ :\ |x|\le t+R\}
\end{equation}
with some $R\ge1$. If
\begin{equation}
\label{cond22}
p>\frac{2n}{n-1}\quad\mbox{and}\quad1<q <\frac{n+1}{n-1},
\end{equation}
then there exists a constant $\e_0=\e_0(f,g,n,p,\mu, \beta, R)>0$
such that $T$ has to satisfy
\begin{equation}
\label{lifespan2}
T(\e)\le C\e^{-(q-1)/\{q+1-n(q-1)\}}
\end{equation}
for $0<\e\leq \e_0$, where $C$ is a positive constant independent of $\e$.
\end{thm}

\begin{rem}
Under the assumption \eqref{cond22}, the lifespan estimate \eqref{lifespan2} is better than \eqref{lifespan1}. For this, we should have
\begin{equation}
\label{condforim}
 \frac{q-1}{q+1-n(q-1)}<
 \frac{2p(q-1)}{2q+2-(n-1)p(q-1)}
\end{equation}
 which is equivalent to
\begin{equation}
\label{condforp}
 p>\frac{2(q+1)}{2(q+1)-(n+1)(q-1)}.
\end{equation}
On the other hand, $q<(n+1)/(n-1)$ is equivalent to
\[
 \frac{2(q+1)}{2(q+1)-(n+1)(q-1)}<\frac{2n}{n-1},
 \]
 which means that assumption \eqref{cond22} guarantees the inequality \eqref{condforim}.
In section 6 we will give the reason why we have to pose the restriction on $p$ in the form
\[
p>\frac{2n}{n-1}
\]
instead of \eqref{condforp}.

\end{rem}

\section{Lower bound of the first functional}
\par\quad
One of the key ingredients to the blow-up result is to get the lower bound of
\[
F_1(t):=\int_{\R^n}u(x, t)\psi(x, t)dx,
\]
where
\begin{equation}
\label{test11}
\psi(x,t):=e^{-t}\phi_1(x),
\quad
\phi_1(x):=
\left\{
\begin{array}{ll}
\d\int_{S^{n-1}}e^{x\cdot\omega}dS_\omega & \mbox{for}\ n\ge2,\\
e^x+e^{-x} & \mbox{for}\ n=1,
\end{array}
\right.
\end{equation}
which was first introduced in Yordanov and Zhang \cite{YZ06}.
Another key point is a multiplier,
\begin{equation}
\label{test1}
m(t):=\exp\left(\mu\frac{(1+t)^{1-\beta}}{1-\beta}\right),
\end{equation}
which is crucial for our proof and was first introduced in \cite{LT18}.
We note that $m(t)$ is bounded as
\[
0<m(0)\leq m(t)\leq 1.
\]
Then we have the following lemma.
\begin{lem}
\label{F1}
Let $u$ be an energy solution of \eqref{nonlinear3} on $[0, T)$. Under the same assumption of Theorem \ref{blowup}, it holds that
\begin{equation}
\label{F1postive}
F_1(t)\ge \frac{m(0)\e}{2}\int_{\R^n}f(x)\phi_1(x)dx\ge0
\quad\mbox{for}\ t\ge0.
\end{equation}
\end{lem}
\par\noindent
{\it Proof.} The proof of Lemma \ref{F1} is almost the same as that of Lemma 3.1 in \cite{LT17}, which is established by neglecting the spatial integral of the nonlinear term
\[
\int_{\R^n}|u(x, t)|^pdx
\]
due to its positivity.
Replacing this quantity by
\[
\int_{\R^n}\left\{|u_t(x, t)|^p+|u(x, t)|\right\}^qdx,
\]
we get the desired proof immediately.
\section{Lower bound of the second functional}
\par\quad
With Lemma \ref{F1} in hand, we may prove a key inequality for
\[
F_2(t):=\int_{\R^n}u_t(x, t)\psi(x, t)dx.
\]
\begin{lem}\label{F2}
Let $u(x, t)$ and $\psi(x, t)$ be as in section 3. Under the same assumption of Theorem \ref{blowup}, it holds that
\begin{equation}
\label{F2postive}
F_2(t)\ge \frac{m(0)\e}{2}\int_{\R^n}g(x)\phi_1(x)dx\ge0
\quad\mbox{for}\ t\ge0.
\end{equation}
\end{lem}
{\it Proof.} Actually Lemma \ref{F2} is a partial result of the proof of Theorem 2.1 in \cite{LT17}. For convenience we rewrite the detail.
By direction calculation we have
\begin{equation}
\label{u_t+u}
\begin{array}{l}
\d\frac{d}{dt}\left[m(t)\int_{\R^n}\left\{u_t(x, t)+u(x, t)\right\}\psi(x, t)dx\right]\\
\d=\frac{\mu}{(1+t)^\beta}m(t)\int_{\R^n}\left\{u_t(x, t)+u(x, t)\right\}\psi(x, t)dx\\
\quad\d+m(t)\frac{d}{dt}\int_{\R^n}\left\{u_t(x, t)+u(x, t)\right\}\psi(x, t)dx.
\end{array}
\end{equation}
Replacing the test function $\phi$ in the definition \eqref{energysol} with $\psi$
and taking derivative to both sides with respect to $t$, we have that
\begin{equation}
\label{derivateofenergysol}
\begin{array}{l}
\d\frac{d}{dt}\int_{\R^n}u_t(x,t)\psi(x,t)dx-\int_{\R^n}u_t(x,t)\psi_t(x,t)dx\\
\d+\int_{\R^n}\nabla u(x,t)\cdot \nabla\psi(x,t)dx
+\frac{\mu}{(1+t)^\beta}\int_{\R^n}u_t(x,t)\psi(x,t)dx\\
\d=\int_{\R^n}|u_t(x,t)|^p\psi(x,t)dx+\int_{\R^n}|u(x,t)|^q\psi(x,t)dx.
\end{array}
\end{equation}
Since for $\psi(x, t)$ we have
\[
\psi_t=-\psi,\quad\psi_{tt}=\Delta\psi=\psi,
\]
then by integration by parts in the first term in the second line
of the last equality yields that
\begin{equation}
\label{IBP}
\begin{array}{l}
\d\frac{d}{dt}\int_{\R^n}\left\{u_t(x,t)+u(x,t)\right\}\psi(x,t)dx\\
\d+\frac{\mu}{(1+t)^\beta}\int_{\R^n}u_t(x,t)\psi(x,t)dx\\
\d=\int_{\R^n}|u_t(x,t)|^p\psi(x,t)dx+\int_{\R^n}|u(x,t)|^q\psi(x,t)dx.
\end{array}
\end{equation}
By combining \eqref{u_t+u} and \eqref{IBP} we have
\begin{equation}
\label{derivateu_t+u}
\begin{array}{l}
\d\frac{d}{dt}\left[m(t)\int_{\R^n}\left\{u_t(x, t)+u(x, t)\right\}\psi(x, t)dx\right]\\
\d=m(t)\int_{\R^n}|u_t(x,t)|^p\psi (x,t)dx+m(t)\int_{\R^n}|u(x,t)|^q\psi (x,t)dx\\
\d+\frac{\mu}{(1+t)^\beta}m(t)F_1(t)
\end{array}
\end{equation}
for $t\ge0$.
Then (\ref{derivateu_t+u}) and the positivity of $F_1$ by Lemma \ref{F1} yield
\begin{equation}
\label{integrationu_t+u}
\begin{array}{l}
\d m(t)\int_{\R^n}\left\{u_t(x, t)+u(x, t)\right\}\psi(x, t)dx\\
\d\geq m(0)\e\int_{\R^n}\{f(x)+g(x)\}\phi_1(x)dx\\
\d\quad
+\int_0^tds\int_{\R^n}m(s)|u_t(x, s)|^p\psi(x, s) dx\\
\d\quad+\int_0^tds\int_{\R^n}m(s)|u(x, s)|^q\psi(x, s) dx.
\end{array}
\end{equation}
\par
On the other hand, noting that
\[
\frac{m'(t)}{m(t)}=\frac{\mu}{(1+t)^\beta},
\]
then \eqref{derivateofenergysol} implies that
\[
\begin{array}{l}
\d\frac{d}{dt}\int_{\R^n}u_t(x,t)\psi(x,t)dx
+\frac{m'(t)}{m(t)}\int_{\R^n}u_t(x,t)\psi(x,t)dx\\
\d+\int_{\R^n}\left\{u_t(x,t)-u(x,t)\right\}\psi(x,t)dx\\
\d=\int_{\R^n}|u_t(x,t)|^p\psi(x,t)dx+\int_{\R^n}|u(x,t)|^q\psi(x,t)dx.
\end{array}
\]
Multiplying the above equality by $m(t)$, we get
\begin{equation}
\label{derivateandintegration}
\begin{array}{l}
\d\frac{d}{dt}\left[m(t)\int_{\R^n}u_t(x,t)\psi(x,t)dx\right]\\
\d+m(t)\int_{\R^n}\left\{u_t(x,t)-u(x,t)\right\}\psi(x,t)dx\\
\d=m(t)\int_{\R^n}|u_t(x,t)|^p\psi(x,t)dx+m(t)\int_{\R^n}|u(x,t)|^q\psi(x,t)dx.
\end{array}
\end{equation}
Adding \eqref{integrationu_t+u} and \eqref{derivateandintegration} together,
we obtain that
\begin{equation}
\label{addtwo}
\begin{array}{l}
\d\frac{d}{dt}\left[m(t)\int_{\R^n}u_t(x,t)\psi(x,t)dx\right]
+2m(t)\int_{\R^n}u_t(x,t)\psi(x,t)dx\\
\d\geq m(0)\e\int_{\R^n}\{f(x)+g(x)\}\phi_1(x)dx
+m(t)\int_{\R^n}|u_t(x,t)|^p\psi(x,t)dx\\
\d\quad+m(t)\int_{\R^n}|u(x,t)|^q\psi(x,t)dx\\
\d\quad+\int_0^tm(s)ds\int_{\R^n}|u_t(x,s)|^p\psi(x,s)dx\\
\d\quad+\int_0^tds\int_{\R^n}m(s)|u(x, s)|^q\psi(x, s) dx.\\
\end{array}
\end{equation}
Setting
\begin{equation}
\label{G}
\begin{array}{ll}
G(t):=
&
\d m(t)\int_{\R^n}u_t(x,t)\psi(x,t)dx
-\frac{m(0)\e}{2}\int_{\R^n}g(x)\phi_1(x)dx
\\
&
\d -\frac{1}{2}\int_0^tm(s)ds\int_{\R^n}|u_t(x,s)|^p\psi(x,s)dx,
\end{array}
\end{equation}
then we have
\[
G(0)=\frac{m(0)\e}{2}\int_{\R^n}g(x)\phi_1(x)dx>0.
\]
It is easy to get from \eqref{addtwo} that
\[
\begin{aligned}
&G'(t)+2G(t)\\
\geq &\frac{m(t)}{2}\int_{\R^n}|u_t(x,t)|^p\psi(x,t)dx
+m(0)\e\int_{\R^n}\phi_1(x)f(x)dx\\
\geq &0
\end{aligned}
\]
which implies
\[
G(t)\geq e^{-2t}G(0)>0 \quad\mbox{for}\ t\ge0.
\]
Hence, by the definition \eqref{G}, it holds that
\begin{equation}
\label{inequalityforG(t)}
\begin{array}{l}
\d m(t)\int_{\R^n}u_t(x,t)\psi(x,t)dx\\
\d \ge\frac{1}{2}\int_0^tm(s)ds\int_{\R^n}|u_t(x,s)|^p\psi(x,s)dx\\
\d\quad+\frac{m(0)\e}{2}\int_{\R^n}g(x)\phi_1(x)dx,
\end{array}
\end{equation}
which implies that
\[
\begin{aligned}
\int_{\R^n}u_t(x,t)\psi(x,t)dx&\geq \frac{m(0)\e}{2m(t)}\int_{\R^n}g(x)\phi_1(x)dx\\
&\geq \frac{m(0)\e}{2}\int_{\R^n}g(x)\phi_1(x)dx,
\end{aligned}
\]
which is exactly the desired inequality in Lemma \ref{F2}.
\section{Iteration argument}
\par\quad
As mentioned in the introduction, we can't establish the blow-up result and lifespan estimate by using Kato's lemma, instead of which we will
use an iteration argument, following the idea in \cite{LT18}. Set
\[
F_0(t):=\int_{\R^n}u(x,t)dx.
\]
Choosing the test function $\phi=\phi(x,s)$ in \eqref{energysol} to satisfy
$\phi\equiv 1$ in $\{(x,s)\in \R^n\times[0,t]:|x|\le s+R\}$, we get
\[
\begin{array}{l}
\d\int_{\R^n}u_t(x,t)dx-\int_{\R^n}u_t(x,0)dx+\int_0^tds\int_{\R^n}\frac{\mu u_t(x,s)}{(1+s)^\beta}dx\\
=\d\int_0^tds\int_{\R^n}|u_t(x,s)|^pdx+\int_0^tds\int_{\R^n}|u(x,s)|^qdx,
\end{array}
\]
which implies that by taking derivative with respect to $t$ on the both sides
\begin{equation}
\nonumber\\
F_0''(t)+\frac{\mu}{(1+t)^\beta}F_0'(t)
=\int_{\R^n}|u_t(x,t)|^pdx+\int_{\R^n}|u(x,s)|^qdx.
\end{equation}
Multiplying with $m(t)$ on the both sides yields
\begin{equation}
\label{F'0eq}
\left\{m(t)F'_0(t)\right\}'
=m(t)\int_{\R^n}|u_t(x,t)|^pdx+m(t)\int_{\R^n}|u(x,t)|^qdx,
\end{equation}
which means that
\begin{equation}
\label{F'0ineq}
F'_0(t)\geq m(0)\int_0^tds\int_{\R^n}|u_t(x,s)|^pdx+m(0)\int_0^tds\int_{\R^n}|u(x,s)|^qdx.
\end{equation}
\begin{lem}[Inequality (2.5) of Yordanov and Zhang \cite{YZ06}]
\label{lem1}
There exists a constant $C_1=C_1(n,p,R)>0$ such that
\begin{equation}
\label{psi}
\int_{|x|\leq t+R}\left[\psi(x,t)\right]^{p/(p-1)}dx
\leq C_1(1+t)^{(n-1)\{1-p/(2(p-1))\}}
\quad\mbox{for}\ t\ge0.
\end{equation}
\end{lem}
By H\"{o}lder's inequality, \eqref{psi} and \eqref{F2}, we may estimate the nonlinear term
\[
\begin{aligned}
\int_{\R^n}|u_t(x,t)|^pdx&\geq F_2^p(t)\left(\int_{|x|\leq t+R}\left[\psi(x,t)\right]^{p/(p-1)}dx\right)^{-(p-1)}\\
&\geq C_2\e^p(1+t)^{-(n-1)(p-2)/2},\\
\end{aligned}
\]
where
\[
C_2:=C_1^{1-p}\left(\frac{m(0)}{2}\int_{\R^n}g(x)\phi_1(x)dx\right)^p.
\]
Plugging which into \eqref{F'0ineq} we have
\begin{equation}
\begin{aligned}\label{F0first}
F_0(t)&\geq m(0)C_2\e^p\int_0^t\int_0^s(1+r)^{n-1-(n-1)p/2}drds\\
&\geq m(0)C_2\e^p(1+t)^{-(n-1)p/2}\int_0^t\int_0^sr^{n-1}drds\\
&\geq C_3\e^p(1+t)^{-(n-1)p/2}t^{n+1},
\end{aligned}
\end{equation}
where
\[
C_3:=\frac{m(0)C_2}{n(n+1)}.
\]
By H\"{o}lder's inequality again, it follows from \eqref{F'0ineq} that
\begin{equation}
\label{iteration}
F_0(t)\geq C_4m(0)\int_0^t\int_0^s(1+r)^{-n(q-1)}F_0^q(r)drds
\end{equation}
with some positive constant $C_4$ independent of $\e$.
In this way, we find two key ingredients for our iteration argument.

Assuming that
\begin{equation}
\begin{aligned}\label{F0lower}
F_0(t)\geq A_j(1+t)^{-a_j}t^{b_j}\quad\mbox{for}\ t\ge0
\quad(j=1,2,3\cdots)
\end{aligned}
\end{equation}
with
\begin{equation}\label{coe1}
A_1=C_3\e^p,~a_1=\frac{(n-1)p}{2},~b_1=n+1.
\end{equation}
Plugging \eqref{F0lower} into \eqref{iteration} we have
\[
\begin{aligned}
F_0(t)\geq A_{j+1}(1+t)^{-qa_j-n(q-1)}t^{qb_j+2},
\end{aligned}
\]
where
\begin{equation}\label{coej+1}
A_{j+1}\geq \frac{C_4m(0)A_j^q}{(qb_j+2)^2},\ a_{j+1}=qa_j+n(q-1),\ b_{j+1}=qb_j+2.
\end{equation}
By combining \eqref{coe1} and \eqref{coej+1} we come to
\[
\begin{aligned}
&a_j=q^{j-1}\left(\frac{(n-1)p}{2}+n\right)-n,\\
&b_j=q^{j-1}\left(n+1+\frac{2}{q-1}\right)-\frac{2}{q-1},\\
&A_j\geq \frac{C_5A_{j-1}^q}{q^{2(j-1)}}
\end{aligned}
\]
with
\[
C_5:=\frac{C_4m(0)}{\left(n+1+\frac{2}{q-1}\right)^2}.
\]
Hence we have
\[
\begin{aligned}
&\log A_j\\
\geq &q\log A_{j-1}-2(j-1)\log q+\log C_5\\
\geq &q^2\log A_{j-2}-2\big(q(j-2)+(j-1)\big)\log q+(q+1)\log C_5.
\end{aligned}
\]
Repeating this procedure, we have
\[
\log A_j
\geq q^{j-1}\log A_1-\sum_{k=1}^{j-1}\frac{2k\log q-\log C_5}{q^k},
\]
which yields that
\[
A_j\ge\exp\left\{q^{j-1}\left(\log A_1-S_q(j)\right)\right\},
\]
where
\[
S_q(j):=\sum_{k=1}^{j-1}\frac{2k\log q-\log C_5}{q^k}.
\]
By d'Alembert's criterion we know that $S_q(j)$ converges for $q>1$ as $j\rightarrow\infty$.
And therefore we obtain that
\[
A_j\ge\exp\left\{q^{j-1}\left(\log A_1-S_q(\infty)\right)\right\}.
\]
So if we come back to \eqref{F0lower} we have
\begin{equation}
\begin{aligned}\label{F0infty}
F_0(t)&\geq A_j(1+t)^{-a_j}t^{b_j}\\
&\geq (1+t)^nt^{-2/(q-1)}\exp\left(q^{j-1}J(t)\right),~~~t>0,
\end{aligned}
\end{equation}
where
\[
\begin{aligned}
J(t)=&-\left((n-1)\frac{p}{2}+n\right)\log(1+t)
+\left(n+1+\frac{2}{q-1}\right)\log t\\
&+\log A_1-S_q(\infty).
\end{aligned}
\]
Then for $t\geq 1$, $J(t)$ can be estimated as
\begin{equation}
\begin{aligned}\label{J(t)}
J(t)\geq& -\Big((n-1)\frac{p}{2}+n\Big)\log(2t)+\big(n+1+\frac{2}{q-1}\big)\log t\\
&+\log A_1-S_q(\infty)\\
=&\left(n+1+\frac{2}{q-1}-(n-1)\frac p2-n\right)\log t+\log A_1\\
&-\left((n-1)\frac{p}{2}+n\right)\log 2-S_q(\infty)\\
=&\log \left(t^{1+2/(q-1)-(n-1)p/2}A_1\right)-C_6,\\
\end{aligned}
\end{equation}
where
\[
 C_6:=\Big((n-1)\frac{p}{2}+n\Big)\log 2+S_q(\infty).
 \]

Recall the definition of $A_1$ in \eqref{coe1}, we have that $J(t)>1$ if
\[
t\geq C_7\e^{-2p(q-1)/\{2q+2-(n-1)p(q-1)\}}
\]
with
\[
C_7:=\left(C_3^{-1}e^{1+C_6}\right)^{2(q-1)/\{2q+2-(n-1)p(q-1)\}}.
\]
By \eqref{F0infty}, it is easy to get
\[
F_0(t)\rightarrow \infty\quad\mbox{as}\quad j\rightarrow \infty.
\]
Hence we get the lifespan estimate in Theorem \ref{blowup}.

\begin{rem}
In the last line of \eqref{J(t)}, we should require that
\[
1+\frac{2}{q-1}-\frac{(n-1)p}{2}>0,
\]
which leads to the restriction \eqref{cond2} for $q$ in the case $n\geq 2$.
\end{rem}
\section{Proof of Theorem \ref{improvedlifespan}}
\par\quad
Due to \eqref{F0first}, we roughly get an estimate of the form,
\[
F_0(t)\ge C\e^pt^{n+1-(n-1)p/2}
\]
for large $t$ with some positive constant $C$ independent of $\e$. So if
\[
p>\frac{2n}{n-1},
\]
then we have
\[
n+1-(n-1)p/2<1,
\]
which means that \eqref{F0first} is weaker than the linear growth. And hence it is natural to get a better result if we
have linear growth in the first step in the iteration argument. Actually, due to the assumption of the initial data, we get from \eqref{F'0eq} that
\[
F_0'(t)\geq \frac{m(0)}{m(t)}F'_0(0)\geq \left(m(0)\int_{\R^n}g(x)dx\right)\e,
\]
which implies that
\begin{equation}
\begin{aligned}\label{lineargrowth}
F_0(t)\geq C_8\e t,~~~t\geq 0,
\end{aligned}
\end{equation}
where
\[
C_8:=m(0)\int_{\R^n}g(x)dx.
\]
Plugging \eqref{lineargrowth} into \eqref{iteration} we obtain
\begin{equation}
\begin{aligned}\label{firststep}
F_0(t)&\geq C_9\e^q\int_0^t\int_0^s(1+r)^{-n(q-1)}r^qdrds\\
&\geq C_{10}\e^q (1+t)^{-n(q-1)}t^{q+2},\\
\end{aligned}
\end{equation}
where
\[
C_9:=C_4m(0)C_8^q
\quad\mbox{and}\quad
C_{10}:=\frac{C_9}{(q+1)(q+2)}.
\]
Then as in section 5, we may assume that
\begin{equation}
\begin{aligned}\label{F0imlower}
F_0(t)\geq \widetilde{A}_j(1+t)^{-\widetilde{a}_j}t^{\widetilde{b}_j}\quad\mbox{for}\ t\ge0
\quad(j=1,2,3\cdots)
\end{aligned}
\end{equation}
with
\begin{equation}\label{coeim1}
\widetilde{A}_1=C_{10}\e^q,~\widetilde{a}_1=n(q-1),~\widetilde{b}_1=q+2.
\end{equation}
Plugging \eqref{F0imlower} into \eqref{iteration} we have
\[
\begin{aligned}
F_0(t)\geq \widetilde{A}_{j+1}(1+t)^{-q\widetilde{a}_j-n(q-1)}t^{q\widetilde{b}_j+2},
\end{aligned}
\]
where
\begin{equation}\nonumber
\widetilde{A}_{j+1}\geq \frac{C_4m(0)\widetilde{A}_j^q}{(q\widetilde{b}_j+2)^2},\ \widetilde{a}_{j+1}=q\widetilde{a}_j+n(q-1),\ \widetilde{b}_{j+1}=q\widetilde{b}_j+2,
\end{equation}
from which we get that
\begin{equation}
\label{coeim+1}
\left\{
\begin{array}{ll}
\widetilde{a}_j=nq^j-n,\\
\widetilde{b}_j=q^{j-1}\left\{q+2+2/(q-1)\right\}-2/(q-1),\\
\widetilde{A}_j\geq C_{11}\widetilde{A}_{j-1}^q/q^{2(j-1)}
\end{array}
\right.
\end{equation}
with
\[
C_{11}:=\frac{C_4m(0)}{\left\{q+2+2/(q-1)\right\}^2}.
\]
In the same way as in section 5, we conclude that
\[
\widetilde{A}_j\geq \exp\left\{q^{j-1}\left(\log \widetilde{A}_1-\widetilde{S}_q(\infty)\right)\right\}
\]
with
\[
\widetilde{S}_q(\infty):=\lim_{j\rightarrow \infty}\widetilde{S}_q(j):=\lim_{j\rightarrow \infty}\sum_{k=1}^{j-1}\frac{2k\log q-\log C_{11}}{q^k},
\]
and
\begin{equation}
\label{F0im}
F_0(t)\geq (1+t)^nt^{-2/(q-1)}\exp\left(q^{j-1}\widetilde{J}(t)\right)
\end{equation}
with
\[
\widetilde{J}(t)=-nq\log(1+t)+\left(q+2+\frac{2}{q-1}\right)\log t+\log\widetilde{A}_1-\widetilde{S}_q(\infty).
\]
Therefore, if $t\geq 1$, we come to
\[
\begin{aligned}
\widetilde{J}(t)&\geq -nq\log(2t)+\left(q+2+\frac{2}{q-1}\right)\log t+\log\widetilde{A}_1-\widetilde{S}_q(\infty)\\
&=\left(q+2+\frac{2}{q-1}-nq\right)\log t+\log\widetilde{A}_1-\widetilde{S}_q(\infty)-nq\log 2\\
&=\log\left(t^{q+2+2/(q-1)-nq}\widetilde{A}_1\right)-C_{12}
\end{aligned}
\]
with
\[
C_{12}:=\widetilde{S}_q(\infty)+nq\log 2.
\]
If
\[
t\geq C_{13}\e^{-(q-1)/\{q+1-n(q-1)\}},
\]
where
\[
C_{13}:=\left(\frac{e^{C_{12}+1}}{C_{10}}\right)^{1/\left\{q+2+2/(q-1)-nq\right\}},
\]
then we have $\widetilde{J}(t)\geq 1$, which will lead to by \eqref{F0im}
\[
F_0(t)\rightarrow \infty\quad\mbox{as}\quad j\rightarrow \infty,
\]
and we finish the proof of Theorem \ref{improvedlifespan}.

\section*{Acknowledgment}
\par\quad
The first author is partially supported by Zhejiang Province
Science Foundation(LY18A010008), NSFC(11501273, 11726612, 11771359,
11771194), Chinese Postdoctoral Science Foundation(2017M620128), the Scientific Research Foundation of the First-Class Discipline of Zhejiang Province
(B)(201601).The second author is partially supported by the Grant-in-Aid for Scientific Research(C)
(No.15K04964),
Japan Society for the Promotion of Science,
and Special Research Expenses in FY2017, General Topics(No.B21),
Future University Hakodate.


\bibliographystyle{plain}

\end{document}